\input amstex
\documentstyle{amsppt}

\input epsf.sty

\magnification=1200
\vsize=7.25in
\hsize=5.25in
\hfuzz 2pt

\define\BbbY{\Bbb Y}
\define\vk{\varkappa}
\define\BbbC{\Bbb C}
\define\PP{\Cal P}
\define\Tab{\operatorname{Tab}}
\define\Cyl{\operatorname{Cyl}}
\define\HH{\Cal H}

\TagsOnRight
\NoRunningHeads

\topmatter

\title Random Young tableaux and combinatorial identities
\endtitle

\author Grigori Olshanski$^*$ and
Amitai Regev$^{**}$
\endauthor

\thanks $^{\phantom{*}*}$ Supported by the Meyerhoff Foundation. \endgraf
$^{**}$ Partially supported by ISF
Grant 6629 and by Minerva Grant No. 8441
\endthanks

\abstract We derive new combinatorial identities which may be viewed
as multivariate analogs of summation formulas for hypergeometric
series. As in the previous paper \cite{Re}, we start with probability
distributions on the space of the infinite 
Young tableaux. Then we calculate the probability that the entry of
a random tableau at a given box equals $n=1,2,\dots$. Summing these
probabilities over $n$ and equating the result to 1 we get a
nontrivial identity. Our choice for the initial distributions is
motivated by the recent work on harmonic analysis on the infinite
symmetric group and related topics.
\endabstract

\endtopmatter

\head \S0. Introduction \endhead

Let $\Tab$ be the set of all infinite standard Young tableaux
$T=(T(i,j))$. Given a probability measure $M$ on $\Tab$, we may speak
about the {\it random\/} infinite tableau $T$. Let $\PP_M(T(i,j)=n)$ denote
the probability that $T$ has the entry $n$ at the box $(i,j)$.
Fix a box $(i,j)$ such that the shape of $T$ contains $(i,j)$ almost
surely. Then 
$$
\sum _ {n \ge 0} \PP_M(T(i,j)=n)=1.\tag0.1
$$
In \cite{Re}, it was shown that by specializing $M$ one can get from (0.1)
many nontrivial identities. These identities look as summation
formulas for multivariate series of hypergeometric type. 

For instance, one of the identities is as follows 
(see \cite{Re, (4.2.$1'$)}): 
$$
\gathered
\sum\Sb p_1>\dots>p_k\ge1\\ q_1>\dots>q_l\ge1\endSb
\frac{(|p|+|q|+\tfrac12[k+l-(k+l)^2])!\,V^2(p)\,V^2(q)}
{\prod\limits_{1\le r\le k}(p_r!)^2\,
\prod\limits_{1\le s\le l}(q_s!)^2\,
\prod\limits_{1\le r\le k}\prod\limits_{1\le s\le l}
(p_r+q_s+1)^2}\\
\times\left(\prod_{r=1}^k\frac{p_r}{p_r+1}\right)
\left(\prod_{s=1}^l\frac{q_s}{q_s+1}\right)=1,
\endgathered \tag0.2
$$
where $k,l=0,1,\dots$ are arbitrary, $k+l\ge1$, and 
$$
\gather
|p|=\sum_{r=1}^k p_r\,, \qquad 
|q|=\sum_{s=1}^l q_s\,, \\
V^2(p)=\prod_{1\le i<j\le k}(p_i-p_j)^2, \qquad
V^2(q)=\prod_{1\le i<j\le l}(q_i-q_j)^2.
\endgather
$$
The identity (0.2) arises from the so--called {\it Plancherel
measure,\/} the fixed box being $(k+1,l+1)$.  
{}From a different point of view, the identity (0.2) (written in an
equivalent form) is discussed in \cite{MMW}. \footnote{See Remark 2.8
for a comment to the approach of \cite{MMW}.}

In the present paper, which is a continuation of \cite{Re}, we derive
new identities of the form (0.1). The results are as follows.

\proclaim{Theorem 0.1} Let $k, l=0,1,2,\dots$ be arbitrary,
$k+l\ge1$, and $t >0$ be a parameter. Then 
$$
\multline
\frac{1}{k!} \sum_{r_1,\dots,r_l,s_1, \dots,s_k \geq 0}
\frac{(s_1+ \dots +s_k+r_1+2r_2+ \dots +lr_l+kl+k+l)!}
{(s_1+l+1) \dots (s_k+l+1)1^{r_1}2^{r_2} \dots l^{r_l}r_1! \dots r_l!} \\
\times \frac{t^{k+r_1+ \dots +r_l+1}}
{(t)_{s_1+ \dots +s_k+r_1+2r_2+ \dots +lr_l+kl+k+l+1}}
=1, 
\endmultline \tag 0.3
$$
where $(x)_n=x(x+1) \dots(x+n-1)$ is the Pochhammer symbol.
\endproclaim

See Theorem 3.3.2 below. We present two proofs of (0.3). One of them
follows our general scheme while another is a direct argument, which
is largely due to S.~Milne \cite{Mi}. 
Of the general identities
derived here and in \cite{Re}, so far this is the only case where an
independent direct proof was found.

\proclaim{Theorem 0.2} Let $k=1,2,\dots$, and let $\theta>0$ and $z\in
\BbbC$ be parameters. Then 
$$
\multline
\sum_{\mu _1 \ge \dots \ge \mu _{k+1} =1}
\frac{|\mu |! \cdot 
\prod_{1 \le i < j \le k}((j-i)\theta +\mu _i-\mu _j)}
{\theta \cdot \prod_{1 \le i \le k} (\mu _i-\mu _{i+1})! \cdot
\prod_{1 \le i < j \le k+1} 
((j-i)\theta +\mu _i-\mu _{j-1})_{\mu _{j-1}-\mu _{j+1}+1}} \\
\times \frac{(z-k\theta )(\bar z-k\theta ) \cdot 
\prod_{1 \le i \le k}
\left[ (z-(i-1)\theta )_{\mu _i}
(\bar z-(i-1)\theta )_{\mu _i} \right]}
{(\theta ^{-1}z\bar z)_{|\mu |+1}} =1\ ,
\endmultline \tag 0.4
$$
where $| \mu |=\mu _1+ \dots + \mu _k\ .$
\endproclaim

See Theorem 2.5.1. For instance, in the simplest case $k=1$ we get
the following special summation formula for the generalized
hypergeometric series of type (3,2): 
$$
{}_3F_2(z+1,\bar z+1, 2; \theta +2, \theta^{-1}z \bar z+2;1)=
\frac{(\theta +1)(z \bar z+\theta )}{(z-\theta )(\bar z- \theta )}\,.
\tag 0.5
$$
This can be derived from a certain known formula, see (2.6.1) below. 

Specializing $\theta =1$ and letting $z \to \infty$, we get from (0.4)
the following identity:
$$
\sum_{\mu_1 \ge \dots \ge \mu_{k+1}=1}
\frac{(\mu_1+ \dots +\mu _k)! 
\prod_{1 \le i < j \le k}(\mu_i-\mu_j+j-i)}
{\prod_{i=1}^k(\mu_i - \mu_{i+1})! \cdot
\prod_{1 \le i < j \le k+1}
(\mu_i-\mu_{j-1}+j-i)_{\mu_{j-1}-\mu_{j+1}+1}}\,,
\tag 0.6
$$
which can be transformed (see Proposition 2.7.1) to
$$
\sum_{p_1 > \dots > p_k \ge 1}
\frac{(p_1+ \dots + p_k -k(k-1)/2)! \cdot
\prod_{1 \le i < j \le k}(p_i-p_j)^2}
{ \prod_{i=1}^k[(p_i-1)!(p_i+1)!] } =1\,. \tag 0.7
$$
The last identity is a particular case of (0.2) (it corresponds to
$l=0$). 

Let us explain now the origin of the measures $M$ that lead to the
identities (0.3) and (0.4).

Infinite Young tableaux can be identified with infinite paths in a
graph $\BbbY$, called the Young graph \cite{VK} (the vertices of
$\BbbY$ are arbitrary Young diagrams). Let $\Cal T$ stand for the
space of all infinite paths in $\Bbb Y$. Vershik and Kerov introduced
in \cite{VK} the concept of a {\it central\/} measure on the space
$\Cal T$. Their definition was inspired by the theory of characters of
the infinite symmetric group $S(\infty)$. However, it makes sense for
certain other graphs as well.  

The probability measures $M$ that were considered in \cite{Re}
are related to indecomposable characters of $S(\infty)$
\cite{KV1, KV2, T, VK} and to some decomposable characters studied in
harmonic analysis on $S(\infty)$ \cite{B2, BO1, BO2, BO3, BO4, KOV, Ro}.
Also considered in \cite{Re} were measures related to {\it
projective\/} characters of $S(\infty)$ \cite{B1, I, N}. 

In the present paper, we deal with a modified definition of central
measures, that of ``$\theta$--central measures''. Here $\theta\ge0$ is
an additional parameter. This definition, due to Kerov \cite{Ke2,
Ke3, Ke5}, is related to an additional structure on $\BbbY$,
namely, to certain formal edge multiplicities depending on the
parameter $\theta$. The formal edge multiplicities in question come
from the Pieri rule for the Jack symmetric functions. \footnote{About
the Jack symmetric functions, see \cite{Ma, S}.} When $\theta = 1$,
all these multiplicities are 
equal to 1, which corresponds to the Schur functions and to the Young
graph. Thus, the class of $\theta$--central measures includes that of
central measures as a particular case. 

Another important particular case is that of $\theta=0$, when the
Jack symmetric functions degenerate to the monomial symmetric
functions. In that case all the edge multiplicities are natural numbers
and we get the so--called Kingman graph. The definition of the
Kingman graph, also due to Kerov (see \cite{Ke1}), was initially
inspired by Kingman's concept of partition structures, see \cite{Ki1,
Ki2}.

In Theorem 0.1, we are dealing with the Kingman graph ($\theta=0$).
The identity in question is related to a remarkable family of
``0--central measures'' depending on a parameter $t>0$. 
These measures come from the well--known Poisson--Dirichlet
distributions studied by many authors, see,
e.g., \cite{Ki3}. The fixed box is $(k+1,l+1)$. Note that one could
write down a generalization of (0.3), which is related to the
so--called two--parameter Poisson--Dirichlet distribution, see
\cite{Ke4, Pi, PY}.

In Theorem 0.2, we are dealing with general Jack edge multiplicities
($\theta>0$). The corresponding ``$\theta$--central measures'' depend
on the complex parameter $z$. These measures are studied in
\cite{BO3, Ke5}. As the fixed box, we take $(k+1,1)$. We did
not consider an arbitrary box only to simplify the presentation of the
identity. 

Our formalism has an evident extension to identities of the form
$$
\sum_{n_1, \dots,n_k}\PP_M(T(i_1,j_1)=n_1, \dots , T(i_k,j_k)=n_k) = 1,
\tag 0.8
$$
which correspond to several boxes. Some examples are mentioned in
\cite{Re, 5.3} and in \cite{MMW}. For simplicity, we do not deal with
this case here. 

\subhead Acknowledgment \endsubhead
The first named author is grateful to the 
Weizmann Institute and the Meyerhoff Foundation for their hospitality
and support during the preparation of this paper.

\head \S1 The general formalism \endhead

\subhead 1.1. The Young graph \endsubhead
Let $\BbbY_n$ be the set of Young diagrams with $n$ boxes and let
$\BbbY=\BbbY_0\sqcup\BbbY_1\sqcup\BbbY_2\sqcup\dots$ be the set of
all Young diagrams. We agree that $\BbbY_0$ consists of a single
element --- the empty diagram $\varnothing$.

For a diagram $\lambda\in\BbbY$ we denote by $|\lambda|$ the number
of boxes in $\lambda$. Given $\mu,\lambda\in\BbbY$, we write $\mu\nearrow$ or
$\lambda\searrow$ if $\mu\subset\lambda$ and $|\lambda|=|\mu|+1$. This
condition means that $\lambda$ is obtained from $\mu$ by adding a
single box. The {\it Young graph\/} is the graph whose
vertices are arbitrary diagrams $\lambda\in\BbbY$ and whose edges
are arbitrary couples $\mu,\lambda\in\BbbY$ such that
$\mu\nearrow\lambda$. By abuse of notation we denote the Young graph
again by the symbol $\BbbY$.

\subhead 1.2. Finite tableaux and paths \endsubhead
Recall that a {\it standard tableau\/} of a given shape
$\lambda\in\BbbY_n$ is defined by labelling the boxes of $\lambda$
with the numbers $1,\dots,n$ in such a way that the labels strictly
increase from left to right along each row and down each column of
$\lambda$. Given a standard tableau $T$ we denote by $T(i,j)$ the label it
assigns to a box $(i,j)\in\lambda$. The set of all standard tableaux
of shape $\lambda$ will be denoted by $\Tab(\lambda)$.

We will identify a standard tableau $T\in \Tab(\lambda)$ with a {\it
path\/}
$\tau=(\varnothing\nearrow\tau^1\nearrow\dots\nearrow\tau^n=\lambda)$.
Here, for any
$k=1,\dots,n$, the diagram $\tau^k$ consists of those boxes
$(i,j)\in\lambda$ for which $T(i,j)\le k$. The correspondence
$T\mapsto\tau$ is a bijection between $\Tab(\lambda)$ and the set of all
paths in the Young graph starting at $\varnothing$ and ending at
$\lambda$.

\subhead 1.3. Infinite tableaux and paths \endsubhead
By an {\it infinite path\/} in the graph $\BbbY$ we mean an
infinite sequence of diagrams
$\tau=(\varnothing\nearrow\tau^1\nearrow\tau^2\nearrow\dots)\;$. We may
view $\tau$ as an {\it infinite standard tableau\/} $T$ whose shape
$D(T)$ is the {\it infinite Young diagram\/} $\bigcup\limits_{n\ge0}\tau^n$.

Let
$$
\Tab_n=\bigcup_{\lambda:\,|\lambda|=n}\Tab(\lambda)
$$
and let $\Tab$ denote the set of all infinite standard tableaux.
For any $n\ge1$ we define the projection $\Tab_n\to \Tab_{n-1}$ as
removing from a given tableau $T\in \Tab_n$ the box labelled by $n$.
In terms of paths $\tau$, this means removing the last diagram $\tau^n$.
Using these projections we can identify $\Tab$ with the projective
limit $\varprojlim \Tab_n$.

We equip $\Tab$ with the topology of a projective
limit of finite sets. In this topology, $\Tab$ is a compact
metrizable totally disconnected topological space. 

Given $\tau\in\Tab_n$, we denote by $\Cyl(\tau)\subset\Tab$ the
pull--back of $\{\tau\}$ under the natural projection
$\Tab\to\Tab_n$. This is a {\it cylindrical subset.\/} Such subsets
are both open and closed, and they form a base of the topology of
$\Tab$. 

\subhead 1.4. Multiplicity function\endsubhead
Assume we are given a strictly positive function $\varkappa(\mu,
\lambda)$ defined on the set of edges of $\BbbY$. Such a function is
called a {\it multiplicity function,\/} and its value at $(\mu,\lambda)$
is called the {\it formal multiplicity\/} of the edge $(\mu,\nu)$. Given
$\tau\in \Tab_n$, we set
$$
\varkappa(\tau)=\varkappa(\varnothing,\tau^1)\varkappa(\tau^1,\tau^2)
\dots \varkappa(\tau^{n-1},\tau^n)\ ,
$$
where $\varnothing\nearrow \tau^1\nearrow \dots 
\nearrow\tau^n$ are the vertices of
$\tau$. We agree that the value of $\varkappa$ at the single element of
$\Tab_0$ is equal to 1.

For a diagram $\lambda$ we set
$$
\dim_{\varkappa}(\lambda)=\sum_{\tau\in \Tab(\lambda)} \varkappa(\tau)\ .
\tag 1.4.1
$$
We call $\dim_{\varkappa}(\cdot)$ the $\varkappa$-{\it dimension function.\/}
If $\varkappa(\cdot,\cdot)\equiv 1$ then $\dim_{\varkappa}(\lambda)$ coincides
with $\vert \Tab(\lambda)\vert$, i.e. it turns into the conventional
combinatorial dimension.\medskip

The $\varkappa$-dimension function satisfies the recurrence relation
$$
\dim_{\varkappa}(\lambda)=\sum_{\mu:\mu\nearrow\lambda}\dim_{\varkappa}
(\mu)\varkappa(\mu,\lambda)\ .\tag 1.4.2
$$
Together with the initial condition $\dim_{\varkappa}(\varnothing)=1$ this
defines the function $\dim_{\varkappa}(\cdot)$ uniquely.

\subhead 1.5. $\varkappa$-central measures\endsubhead
A probability measure $M$ on $\Tab$ is called $\varkappa$-{\it
central\/} if for any diagram $\lambda$, the masses of the cylindrical
sets $\Cyl(\tau)$ with $\tau\in \Tab(\lambda)$ are proportional to the
numbers $\vk(\tau)$, i.e.,
$$
M(\Cyl(\tau))=\vk(\tau)\varphi(\lambda)\ ,\tag 1.5.1
$$
where $\varphi(\lambda)$ is a certain (positive real valued) function on
the Young diagrams. It is easy to see that $\varphi$ is nonnegative,
$\varphi(\varnothing)=1$ (since $\vk(\varnothing)=1=M(\Cyl(\varnothing))$) and for all $\mu,$
$$
\varphi(\mu)=\sum_{\lambda:\lambda\searrow\mu} \vk(\mu,\lambda)
\varphi(\lambda). \tag 1.5.2
$$
Indeed, to see this, pick $\tau\in\Tab(\mu)$ and multiply both sides
by $\vk(\tau)$. Then in the left--hand side we get $M(\Cyl(\tau))$,
and in the right--hand side we get $\sum M(\Cyl(\tau'))$, where the
summation is taken over all paths $\tau'$ obtained from $\tau$ by
adding one extra edge. Since $\Cyl(\tau)$ is the disjoint union of
the sets $\Cyl(\tau')$, we get $M(\Cyl(\tau))=\sum M(\Cyl(\tau'))$,
which is equivalent to (1.5.2).

The relation (1.5.2) is called the $\vk$-{\it harmonicity condition.\/}
The correspondence $M\mapsto\varphi$ is a bijection between the
$\vk$-{\it central\/} probability measures on $\Tab$ and the nonnegative,
$\vk$-harmonic functions on $\BbbY$, normalized at $\varnothing\in\BbbY$.
\smallskip
In the particular case $\vk(\cdot,\cdot)\equiv 1$, $\vk$-central measures
and $\vk$-harmonic functions turn into central measures and harmonic
functions, respectively, as defined in \cite{VK}.

\subhead 1.6. Transition and cotransition probabilities\endsubhead
Let $M$ be a $\vk$-central measure on $\Tab$. We view $(\Tab,M)$ as a
probability space, which makes it possible to speak about the random
element of $\Tab$, i.e., the {\it random infinite tableau\/} or,
equivalently, the {\it random infinite path.}

Given $\lambda\in\BbbY$, the probability that the random path passes
through $\lambda$ equals $$\dim_{\vk}\lambda\cdot\varphi(\lambda).
$$
\noindent Let $(\mu,\lambda)$ be an edge. The probability that the random path
passes through $\mu,$ conditional that it passes through $\lambda$, equals
$$
q(\mu,\lambda)=\frac{\dim_{\vk}\mu\cdot\vk(\mu,\lambda)}
{\dim_{\vk}\lambda} \,. \tag1.6.1
$$
Note that $q(\mu,\lambda)$ does not depend on $M$. This is the so-called
{\it cotransition function.\/} Thus, all $\vk$-central measures with a
fixed $\vk$ have one and the same cotransition function.

The {\it transition function\/} $p_M(\mu,\lambda)$ is, by definition, the
probability that the random path passes through $\lambda,$ conditional that
it passes through $\mu$. We have
$$
p_M(\mu,\lambda)=\frac{\vk(\mu,\lambda)\varphi(\lambda)}
{\varphi(\mu)}\ . \tag 1.6.2
$$
This function is well defined, provided that $\varphi(\mu)>0$, and  it
depends on $M$.

\subhead 1.7. Reachable boxes\endsubhead
Given an infinite tableau $T\in \Tab$, we denote by $D(T)$ its shape, which
is an infinite diagram. A box $(i,j)$ is called {\it completely
reachable\/} (with respect to $M$) if it is contained in $D(T)$ for
almost all $T\in \Tab$ (with respect to $M$). Call $(i,j)$ {\it partially
reachable\/} if the same event holds with a nonzero probability. Note
that this terminology differs from that of \cite{Re, section 2.5}.

Let $D(M)$ and $\tilde D(M)$ stand for the set of all completely
reachable and partially reachable boxes, respectively. Then $D(M)\subseteq
\tilde D(M)$, and both sets are infinite diagrams. The set $\tilde D(M)$
of partially reachable boxes is easily described in terms of $\varphi$:
$$
\tilde D(M)= \bigcup_{\lambda:\varphi(\lambda)>0} \lambda\ .
$$
As for the set $D(M)$, which is more interesting for us, its structure
is not so evident. It may well happen that $D(M)\ne \tilde D(M)$. For
instance, take $M=\frac{1}{2}(M_1+M_2)$, where $M_1$ is the
delta-measure concentrated on the infinite path going along the first
row $i=1$, while $M_2$ is defined similarly by replacing the first row
by the first column $j=1$. Both $M_1$ and $M_2$ are $\vk$-central for
any choice of $\vk$. Consequently, $M$ is $\vk$-central, too. Then
$\tilde D(M)$ is the hook formed by the first row and the first column,
while $D(M)$ consists of a single box $(1,1)$.

\subhead 1.8. The identity ---  general form \endsubhead
Given a box $(i,j)$, we denote by $H^{\prime}(i,j)$ the set of finite
diagrams $\mu$ such that $\mu$ does not contain $(i,j)$ and the shape
$\mu\cup (i,j)$ is a diagram. Here we prove

\proclaim{Theorem 1.8.1}  Let $M$ be a $\vk$-central measure, $\varphi$
the corresponding $\vk$-harmonic function and let $(i,j)$ be a completely
reachable box.

Then
$$
\sum_{\mu\in H^{\prime}(i,j)}\dim_{\vk}\mu \cdot \vk(\mu,\mu
\cup (i,j))\cdot \varphi(\mu\cup(i,j))=1\ . \tag 1.8.1
$$
\endproclaim

In fact, Theorem 1.8.1 is a consequence of the slightly more general
Theorem 1.8.3, which is given below.

Let $\Cal T(i,j)$ be the set of the (infinite) paths
$\tau=(\varnothing\nearrow \tau^1
\nearrow \dots)$ such that the diagram of $\tau$ contains $(i,j)$. For
any $\tau\in\Cal T(i,j)$ there exists a unique $n$ such that $(i,j)$ is
contained in $\tau^n$ (hence in $\tau^{n+1},\tau^{n+2},\ldots)$, but
not in $\tau^{n-1}$. Denote this by $\tau (i,j) =n.$
Given a probability measure $M$ on $\Tab,$ we consider the
$M$--probability that a random path $\tau$ satisfies $\tau(i,j)=n,$
and we denote it by $\PP_M(\tau(i,j)=n).$

Let $\Cal T(i,j;n) \subseteq \Tab$ denote the set of
the infinite tableaux (paths) $\tau$ such that $\tau^{n-1}\in H^{\prime}(i,j)$,
$\tau^n=\tau^{n-1}\cup(i,j)$ (i.e such that $\tau(i,j)=n$):
$$
\Cal T(i,j;n)=\{\tau \in \Tab\,\vert\,\tau^{n-1}\in H^{\prime}(i,j)
\;\, \hbox{\rm and} \;\,
\tau^n=\tau^{n-1}\cup (i,j)\}\ .
$$
Clearly, $\PP_M(\tau(i,j)=n) = M(\Cal T(i,j;n)).$

\proclaim{Lemma 1.8.2}
We have
$$
\PP_M(\tau(i,j)=n)=M(\Cal T(i,j;n))=\sum_
{\mu\in H^{\prime}(i,j)\atop \vert\mu\vert =n-1}
\dim_{\vk}\mu \cdot\vk(\mu,\mu\cup(i,j))\cdot\varphi(\mu\cup(i,j))\
.\tag 1.8.2
$$
\endproclaim

\demo{Proof} Given $\mu\in H^{\prime}(i,j)$ with $\vert\mu\vert
=n-1$, let
$$
\Cal T(i,j;\mu)=\{\tau \in \Tab\,\vert\,\tau^{n-1}=\mu\
\hbox{\rm and}
\ \tau^n=\tau^{n-1}\cup
(i,j)\}\, .\tag 1.8.3
$$
Clearly,
$$
\Cal T(i,j;n)=\bigsqcup_{\mu\in H^{\prime}(i,j)\atop \vert\mu\vert =n-1}
\Cal T(i,j;\mu)\ ,
$$
where $\bigsqcup$ indicates a  disjoint union, hence
 $$
M(\Cal T(i,j;n))=\sum_{\mu\in H^{\prime}(i,j)\atop \vert\mu\vert =n-1}
M(\Cal T(i,j;\mu))\ .
$$
Now the set $\Cal T(i,j;\mu)$ is the disjoint union of
$(d_{\mu})$
cylindrical sets $\Cyl(\sigma)$, where $\sigma$ is an arbitrary path
$\varnothing\nearrow\sigma^{\prime}\nearrow \dots 
\nearrow\sigma^{n-1}\nearrow
\sigma^n$ such that $\sigma^{n-1}=\mu$, $\sigma^n=\mu\cup(i,j)$. 
Here $d_{\mu}=\dim \mu$ is the number of standard tableaux of shape $\mu.$

By the
very definition, $M(\Cyl(\sigma))=\prod^n_{i=1}\vk(\sigma^{i-1},\sigma^i)
\cdot\varphi(\sigma^n)$. This can be written as the product of two
expressions
$$
M(\Cyl(\sigma))=\left [\prod^{n-1}_{i=1}\vk(\sigma^{i-1},\sigma
^i)\right ]\cdot\left [\vk(\mu,\mu\cup(i,j))\cdot\varphi(\mu\cup(i,j))\right
]. \tag 1.8.4
$$
The second expression does not depend on $\sigma$, while summing the
first expression over all possible $(d_{\mu})\;\,\sigma$ gives $\dim_{\vk}
(\mu)$. This concludes the proof.\hfill \qed
\enddemo

Recall that $\Cal T(i,j)=\{\tau\in \Tab \;\vert \; (i,j)\in D(\tau)\}.$
Thus $M(\Cal T(i,j))$
is the probability of a random tableau $\tau$ to have $(i,j)\in
D(\tau)$. Denote that probability by $\PP_M(i,j)$.

\proclaim{Theorem 1.8.3} We have
$$
\PP_M(i,j)=\sum_{\mu\in H
^{\prime}(i,j)}\dim_{\vk}\mu\cdot\vk(\mu,\mu\cup(i,j))\cdot\varphi(\mu\cup
(i,j))\,. \tag 1.8.5
$$
\endproclaim

\demo{Proof}
Clearly, $\Cal T(i,j)=\bigsqcup_{n\geq 1}\Cal T(i,j;n)$, a disjoint union,
hence
$$
\PP_M(i,j)=M(\Cal T(i,j))=\sum^{\infty}_{n=1}M(\Cal T(i,j;n))
$$
and by the above lemma,
$$
\PP_M(i,j)=\sum_{\mu\in H
^{\prime}(i,j)}\dim_{\vk}\mu\cdot\vk(\mu,\mu\cup(i,j))\cdot\varphi(\mu\cup
(i,j))\,. \tag 1.8.6
$$
\hfill \qed
\enddemo

If $(i,j)\in D(M)$ (i.e. $(i,j)$ is completely reachable), $\PP_M(i,j)
=1$, and we have

$$
\sum_{\mu\in H^{\prime}(i,j)}\dim_{\vk}\mu\cdot \vk(\mu,\mu\cup(i,j))\cdot
\varphi(\mu\cup(i,j))=1 \ .
$$
This completes the proof of Theorem 1.8.1.

\head \S2 The Jack graph case \endhead

\subhead 2.1. Jack edge multiplicities \endsubhead
Fix a positive parameter $\theta$. Let $P_{\mu}$ denote the Jack symmetric
function with parameter $\theta$ and index $\mu$. Here $\mu$ is a Young
diagram, and the normalization of $P_{\mu}$ is that of \cite{Ma, VI.10}.
Note that Macdonald uses as the parameter $\alpha=\theta^{-1}$.

The simplest case of Pieri's formula for the Jack symmetric functions has
the form
$$
P_{\mu}P_{(1)}=\sum_{\lambda:\lambda\searrow\mu}\vk_{\theta}(\mu,\lambda)
P_{\lambda}\ ,\tag 2.1.1
$$
where $\vk_{\theta}(\mu,\lambda)$ are certain strictly positive numbers.
An explicit
expression for $\vk_{\theta}(\mu,\lambda)$ is as follows (see \cite{Ma,
VI.10, VI.6}). Let $(i,j)$ be the box $\lambda\backslash\mu$. Then
$$
\vk_{\theta}(\mu,\lambda)=\prod^{i-1}_{k=1} \frac{(a(k,j)+(\ell(k,j)+2)
\theta)\cdot (a(k,j)+1+\ell(k,j)\theta)}
{(a(k,j)+(\ell(k,j)+1)\theta)\cdot (a(k,j)+1+(\ell(k,j)+1)\theta)} \tag 2.1.2
$$
where
$$
a(k,j)=\mu_k-j\ ,\quad \ell(k,j)=\mu^{\prime}_j-k
$$
are the arm--length and the leg--length of the box $(k,j)$ in $\mu$.

We will interprete the numbers $\vk_{\theta}(\mu,\lambda)$ as formal
multiplicities of the edges of the Young graph. By the {\it Jack
graph\/} we mean the Young graph together with this additional
structure. Of course, the Jack graph is not a graph in the
conventional sense.

In the particular case $\theta=1$ the Jack function $P_{\mu}$ turns into the
Schur function $S_{\mu}$, and all the edge multiplicities are equal
to 1, so that the Jack graph reduces to the ordinary Young graph.

The dimension function $\dim_{\vk}\mu$ corresponding to $\vk(\cdot,\cdot)
=\vk_{\theta}(\cdot,\cdot)$ will be denoted as $\dim_{\theta}\mu$. The
following formula is a generalization of the classical hook formula:
$$
\dim_{\theta}\mu=\frac{\vert\mu\vert !}{\HH_{\theta}(\mu)}\ ,\tag 2.1.3
$$
where
$$
\HH_{\theta}(\mu)=\prod_{b\in\mu}(a_{\mu}(b)+\theta\ell_{\mu} (b)+1)
=\prod_{b\in\mu}h _{\mu}(b) 
\tag 2.1.4
$$
and where $h _{\mu}(b)=a_{\mu}(b)+\theta\ell_{\mu} (b)+1\,.$
Here $b\in\mu$ means a box $b=(i,j)$ of $\mu$.

A proof of this formula can be obtained from \cite{S, Theorem 5.4} or
\cite{Ma, VI.10}. A different proof is given in \cite{Ke5, Corollary
(6.10)}. We shall need another expression, which is similar to
$\HH_{\theta}(\mu)$ but different from it:
$$
\HH^{\prime}_{\theta}(\mu)=\prod_{b\in\mu}
(a_{\mu}(b)+\theta\ell_{\mu} (\theta)+\theta)=\prod_{b\in\mu}h^{\prime}  
_{\mu}(b)  \tag 2.1.5 
$$
where $h^{\prime}  _{\mu}(b)=a_{\mu}(b)+\theta\ell_{\mu} (\theta)+\theta\,.$
In the particular case $\theta=1$, $\HH_{\theta}(\mu)$ and $\HH^{\prime}_{
\theta}(\mu)$ coincide.

Finally, for a box $b=(i,j)$, we denote
$$
c_{\theta}(b)=(j-1)-(i-1)\theta\ .
$$
This is the ``$\theta$-version'' of the conventional content $c(b)=j-i$
of a box $b=(i,j).$

Alternative expressions for $\HH_{\theta} (\mu)$ and 
$\HH^{\prime}_{\theta} (\mu)$:
$$
\HH_{\theta}(\mu)=\prod_{1\leq i<j\leq \ell(\mu)}\frac{(1+(j-i-1)\theta)_{\mu_i-\mu_j}}
{(1+(j-i)\theta)_{\mu_i-\mu_j}} \cdot \prod^{\ell(\mu)}_{i=1} (1+(\ell
(\mu)-i)\theta)_{\mu_i} \tag 2.1.6
$$
and
$$
\HH^{\prime}_{\theta} (\mu)=\prod_{1\leq i<j\leq \ell(\mu)}\frac{((j-i)\theta)
_{\mu_i-\mu_j}}
{((j-i+1)\theta)_{\mu_i-\mu_j}} \cdot \prod^{\ell(\mu)}_{i=1} ((\ell
(\mu)+1-i)\theta)_{\mu_i}\,.\tag 2.1.7
$$

\subhead 2.2. The $z$-measures and the Plancherel measures\endsubhead
We fix $\theta >0$ and deal with the multiplicity function $\vk(\cdot,
\cdot)=\vk_{\theta}(\cdot,\cdot)$ as defined above.

\proclaim{Theorem 2.2.1} For any $z\in\BbbC$ there exists a
$\vk_{\theta}$-central measure $M_z$ with the corresponding
$\vk_{\theta}$-harmonic function
$$
\varphi_z(\lambda)=\frac{1}{(\theta^{-1}\vert z\vert^2)_n}\cdot
\frac{\prod_{b\in\lambda}\vert z+c_{\theta}(b)\vert^2}{\HH^{\prime}_{\theta}
(\lambda)}\ ,\tag 2.2.1
$$
where $n=\vert\lambda\vert$, and $(x)_n=x(x+1)\dots (x+n-1)$ stands for
the Pochhammer symbol.
\endproclaim

We call $M_z$ the $z$-{\it measure.\/} About the proof
of this claim see \cite{Ke5} and \cite{BO3, \S3}. In the particular
case $\theta=1$ these measures are discussed in \cite{B2, BO1, BO2,
BO4, KOV}.

\proclaim{Theorem 2.2.2} There exists a $\vk_{\theta}$-central measure
$M_{\infty}$ such that the corresponding $\vk_{\theta}$-harmonic function
$\varphi_{\infty}$ is the pointwise limit of $\varphi_z$ as $\vert z\vert
\to\infty$. We have
$$
\varphi_{\infty}(\lambda)=\frac{\theta^n}{\HH^{\prime}_{\theta} (\lambda)}\ ,
\quad n=\vert\lambda\vert\ .\tag 2.2.2
$$
\endproclaim

The measure $M_{\infty}$ is called the {\it Plancherel measure.\/}
See \cite{KV1, KV2, VK} for the case of the Young graph and
\cite{Ke5} for the general case.

\subhead 2.3. Reachability\endsubhead
Here we use the definitions introduced in section 1.7.

\proclaim {Theorem 2.3.1} Let $M=M_{\infty}$ or $M=M_z$, where $z\in\BbbC$,
$z\notin
\Bbb Z+\Bbb Z\theta$. Then all boxes are completely reachable.
\endproclaim

To prove the theorem we need the following two lemmas.\medskip

\proclaim{Lemma 2.3.2} Fix a box $(i,j)$, let $\mu\in H^{\prime}(i,j)$
be arbitrary, and let $\lambda=\mu\cup (i,j)$. Let $p_{\infty}(\cdot ,
\cdot)$ denote the transition probability for the Plancherel measure
$M_{\infty}$. Then there exists $\epsilon=\epsilon (i,j,\theta)>0,$ not
depending on $\mu,$ such that
$$
p_{\infty}(\mu,\lambda)\geq \epsilon\ .
$$
\endproclaim

\demo{Proof} Recall the general formula
$$
p(\mu,\lambda)=\frac{\vk(\mu,\lambda)\varphi(\lambda)}{\varphi(\mu)}\,.
\tag  2.3.1
$$
When $\vk(\mu,\lambda)=\vk_{\theta}(\mu,\lambda)$, it follows from the
explicit expression for $\vk_{\theta}(\mu,\lambda)$ that this is the
product of $i-1$ factors, each of the type
$$
\frac{(a+\ell\theta +2\theta)(a+\ell\theta +1)}{(a+\ell\theta +\theta)
(a+\ell\theta+\theta +1)}\ .
$$
Since $\theta >0$ is fixed and $a+\ell\theta\geq 0$, this is bounded from
below by a positive constant, uniformly on $a+\ell\theta$. It remains to
examine $\varphi(\lambda)/\varphi(\mu)=\varphi_{\infty}(\lambda)/\varphi
_{\infty}(\mu)$. By (2.2.2),
$$
\frac{\varphi_{\infty}(\lambda)}{\varphi_{\infty}(\mu)}=
\frac{\theta\cdot \HH^{\prime}_{\theta} (\mu)}{\HH^{\prime}_{\theta} 
(\lambda)}\ .\tag 2.3.2
$$
Recall that
$$
\HH^{\prime}_{\theta}(\lambda)=\prod_{b\in\lambda}(a_{\lambda}(b)+\theta\ell
_{\lambda}(b)+\theta)
$$
$$
\HH^{\prime}_{\theta} (\mu)=\prod_{b\in\mu}(a_{\mu}(b)+\theta\ell
_{\mu}(b)+\theta)\ ,
$$
where the subscript indicates that the arm-length and the leg-length are
taken in the corresponding diagram.

Since $\lambda=\mu\cup(i,j)$, we have $a_{\lambda}(b)=a_{\mu}(b)$,
$\ell_{\lambda}(b)=\ell_{\mu}(b)$ whenever $b$ does not lie in the same
column or in the same row that $(i,j)$. Note also that for $b=(i,j)\in
\lambda$ we have $a_{\lambda}(b)=\ell_{\lambda}(b)=0$,  so that $a_{\lambda}
(i,j)+\theta\ell_{\lambda}(i,j)+\theta=\theta$, which cancels out with
$\theta$ in the numerator. Thus, we get
$$
\frac{\varphi_{\infty}(\lambda)}{\varphi_{\infty}(\mu)} =
\prod^{i-1}_{k=1} \frac{a_{\mu}(k,j)+\theta\ell_{\mu}(k,j)+\theta}
{a_{\mu}(k,j)+\theta\ell_{\mu} (k,j)+2\theta} \times
\prod^{j-1}_{k=1} \frac{a_{\mu}(i,k)+\theta\ell_{\mu} (i,k)+\theta}
{a_{\mu}(i,k)+\theta\ell_{\mu} (i,k)+\theta+1}\,.\tag 2.3.3
$$
Clearly, this expression is also bounded from below by a positive constant
not depending on $\mu$.\hfill \qed
\enddemo

\proclaim{Lemma 2.3.3} Fix a box $(i,j)$, let $\mu\in H^{\prime}(i,j)$ be
arbritary, and let $\lambda=\mu\cup(i,j)$. Let $p_z(\mu,\lambda)$ denote
the transition probability for the $z$-measure $M_z$. Assume that
$z\notin\Bbb Z+\Bbb Z\theta$. Then there exists $\epsilon=
\epsilon(i,j,\theta,z)>0$ not depending on $\mu,$ such that
$$
p_z(\mu,\lambda)\geq \frac{\epsilon}{n}\ ,\quad n=\vert\lambda\vert\ .
$$
\endproclaim

\demo{Proof} Comparing the formulas for $\varphi_z$ and
$\varphi_{\infty}$ we see that
$$
p_z(\mu,\lambda)=\frac{\vert z+(j-1)-(i-1)\theta\vert^2}{\vert z\vert^2
+\theta(n-1)} \cdot  p_{\infty}(\mu,\lambda) \ .\tag 2.3.4
$$
By the assumptions on $z$, we have $\vert z+(j-1)-(i-1)\theta\vert >0$. So, the
claim follows from Lemma 2.3.2. \hfill \qed
\enddemo

\demo{Proof of Theorem 2.3.1} Given an infinite Young diagram $D$,
let $\Cal T(D)$ denote the set of all
paths $\tau=(\varnothing\nearrow\tau^1\nearrow\dots)$
such that $D(\tau)=D$. We call $D$ {\it proper\/} if $D$ is distinct from
the set of all the boxes. Then the claim of the theorem will follow
once we prove that for any proper $D,$ $M(\Cal T (D))=0$. So, fix a proper
infinite diagram $D$. There exists a box $(i,j)$ adjacent to $D,$ i.e.
$(i,j)$ does not lie in $D$ but $D\cup (i,j)$ is a diagram. Given $\mu\in
H^{\prime}(i,j)$, let $\Cal T(D,\mu)$ denote the subset of those $\tau\in
\Cal T(D)$ which pass through $\mu$, i.e. $\tau^{\vert\mu\vert}=\mu$.
Clearly, $\Cal T(D)=\bigcup_{\mu\in H^{\prime}(i,j)}\Cal T(D,\mu),$
a countable union. Thus,
it suffices to prove that $M(\Cal T(D,\mu))=0$ for any $\mu\in H^{\prime}
(i,j)$. Set $m=\vert\mu\vert$. For any $n\geq m$ let $\Cal T_n(i,j,\mu)$
denote the set of all $\tau\in\Cal T$ such that $\tau^m=\mu$ and $\tau^n
\in H^{\prime}(i,j)$. Clearly, the sets $\Cal T_n(i,j,\mu)$ decrease
as $n\to\infty$, and for $n\ge m,$ $\Cal T_n(i,j,\mu)$ contains
$\Cal T(D,\mu)$. Thus, it suffices to prove that
$$
\lim_{n\to\infty} M(\Cal T_n(i,j,\mu))=0\ .
$$
We have
$$
M(\Cal T_{m+1}(i,j,\mu))=M(\Cal T_m(i,j,\mu))\cdot (1-\PP_M(\mu,\mu\cup
(i,j)))\ .
$$
By Lemmas 2.3.2, 2.3.3, we have an estimate
$$
\PP_M(\mu,\mu\cup(i,j)) \geq \epsilon_{\vert\mu\vert +1}\ ,
$$
where
$$
\epsilon_n=\cases \epsilon\ , &\text{case of $M=M_{\infty}$}\\
\frac{\epsilon}{n}\ , &\text{case of $M=M_z\ .$}\endcases
$$
Therefore
$$
M(\Cal T_{m+1}(i,j,\mu))\leq (1-\epsilon_{\vert\mu\vert +1})M(\Cal T_m
(i,j,\mu))\ .
$$
Repeating the same argument for any $n>m$, we get
$$
M(\Cal T_n(i,j,\mu))\leq (1-\epsilon_n)M(\Cal T_{n-1}(i,j,\mu))\ .
$$
Since the series $\sum\epsilon_n$ is divergent, it follows that $M(
\Cal T_n(i,j,\mu))\to 0$. \hfill \qed
\enddemo

Theorem 2.3.1 justifies the application of the general Theorem 1.8.1 in
Theorems 2.4.1 and 2.5.1 below. The argument of Theorem 2.3.1 is also
useful to check complete reachability for the examples considered in
\cite{Re}.  

\subhead 2.4. The $\theta$-Plancherel identity for $(k+1,1)$
\endsubhead
Here we specialize $(i,j)=(k+1,1),$ and Theorem 1.8.1 to the case of the
Jack graph,
and we deduce \smallskip

\proclaim{Theorem 2.4.1}
Let $\mu=(a)=(a_1, \dots ,a_k), \quad a_1 \ge \dots \ge a_k \ge
1=a_{k+1}$, $|\mu|=|a|=a_1+ \dots + a_k =n-1, \quad (\lambda=(a_1,
\dots, a_k,1), \quad |\lambda|=n).$  As usual, $a_r=0$ if $r>k+1$.

Then the $\theta$--Plancherel identity for the
box $(k+1,1)$ is
$$
\sum_{a_1 \ge \dots \ge a_{k+1} =1}
\frac{|a|! \cdot \theta^{|a|} \cdot \prod _{1 \le i < j \le k}((j-i)\theta
+a_i-a_j)}
{\prod _{1 \le i \le k} (a_i-a_{i+1})! \cdot
\prod _{1 \le i < j \le k+1} ((j-i)\theta +a_i-a_{j-1})_{a_{j-1}-a_{j+1}+1}}
=1.\tag 2.4.1
$$
\endproclaim

\demo{Proof} We apply Theorem 1.8.1 (the complete reachability is
ensured by Theorem 2.3). Let
$\mu=(\mu_1,\ldots,\mu_k)\vdash n-1$, $\mu_k\geq 1$, 
let $\lambda =
(\mu_1,\ldots, \mu_k,1)\vdash n$, and denote (by abuse of notation)
$\mu_{k+1}=1$. Also $\mu_r=0$ if $r\geq k+2$. Recall that
$$
\HH_{\theta}(\mu)=
\prod_{b\in\mu}(a_{\mu}(b)+\theta\ell_{\mu}(b)+1),\quad \dim_{\theta}\mu=
\frac{\vert\mu\vert !}{\HH_{\theta}(\mu)}\,,\quad 
\HH^{\prime}_{\theta}(\lambda)=
\prod_{b\in\lambda}(a_{\lambda}(\ell)+\theta\ell_{\lambda}(b)+\theta)
$$
and $\varphi_{\infty}(\lambda)=\varphi_{\theta}(\lambda)=\frac{\theta^{
\vert\lambda\vert}}{\HH^{\prime}_{\theta}(\mu)}\,$. Thus, by Theorem 1.8.1, for
the box $(k+1,1)$ we deduce the identity
$$
\sum_{\mu_1\geq\dots\geq\mu_k\geq 1}\frac{\vert\mu\vert !}{\HH_{\theta}
(\mu)}\cdot \vk_{\theta}(\mu,\mu \cup (k+1,1))\cdot\frac{\theta^{\vert\mu\vert +1}}
{\HH^{\prime}_{\theta}(\lambda)}=1\ .\tag 2.4.2
$$
To calculate the left hand side of (2.4.2) explicitly, split the boxes
in $\lambda$ --
hence in $\mu$ -- into the following disjoint subsets
$$
\Delta(i,r)=\{(i,j)\,\vert\,\mu_{r+1}+1\leq j\leq\mu_r\}\ ,\quad 1\leq i\leq r
\leq k+1\ .
$$

\medskip
\midinsert
\epsfxsize 8truecm
\vskip 5pt
\centerline{\epsffile{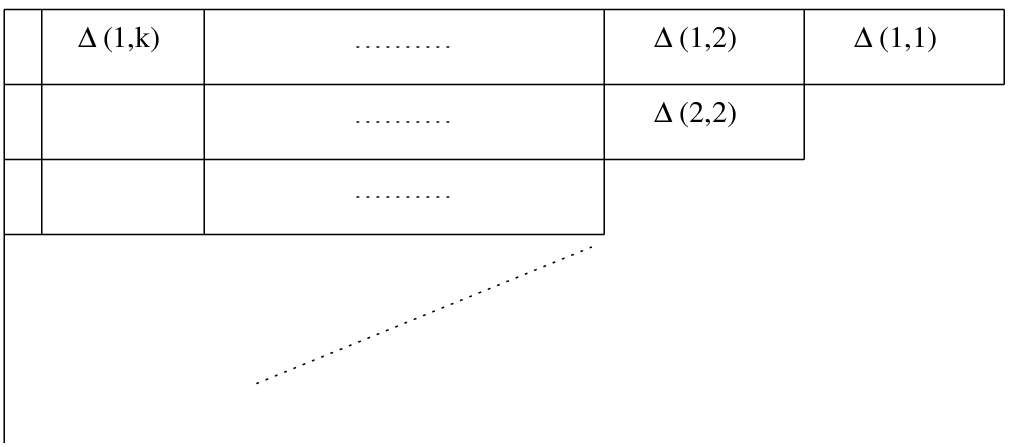}}
\vskip 5pt
\endinsert
\medskip

Thus, for $1\le i \le k,$ the $i$-th row in either $\mu$ or
$\lambda$ is the disjoint union
$\bigsqcup^{k+1}_{r=1}\Delta(i,r)$. Note that $\Delta(i,k+1)=\{(i,1)\}$.
Also,
$
\Delta(k+1,k+1)=\{(k+1,1)\}\ ,\quad (k+1,1)\in\lambda\ $ and
$(k+1,1)\notin\mu .$

Denote
$$
\lambda_*=\bigsqcup_{1\leq i\leq r\leq k-1} \Delta(i,r)\ ,\quad
\mu_*=\bigsqcup_{1\leq i\leq r\leq k-1} \Delta(i,r+1)\ .
$$
Clearly,
$$
\multline
\lambda\backslash\lambda_* = \left (\bigsqcup^k_{i=1}\Delta (i,k)\right )
\sqcup \left (\bigsqcup^{k+1}_{i=1}\Delta(i,k+1)\right ) \\
= \Delta(k+1,k+1) \sqcup \left (\bigsqcup^k_{i=1}(\Delta (i,k)
\sqcup \Delta(i,k+1))\right )\\
= \{(k+1,1)\}\cup \{(i,j)\vert 1\leq i\leq  k,1\leq j\leq \mu_k\},
\endmultline \tag 2.4.3
$$
and
$$
\mu\backslash\mu_* = \left (\bigsqcup^k_{i=1}\Delta(i,i)\right )\sqcup
\{(i,1)\vert 1\leq i\leq k\}\ ,\tag 2.4.4
$$
disjoint unions.

\medskip
\midinsert
\epsfxsize 11truecm
\vskip 5pt
\centerline{\epsffile{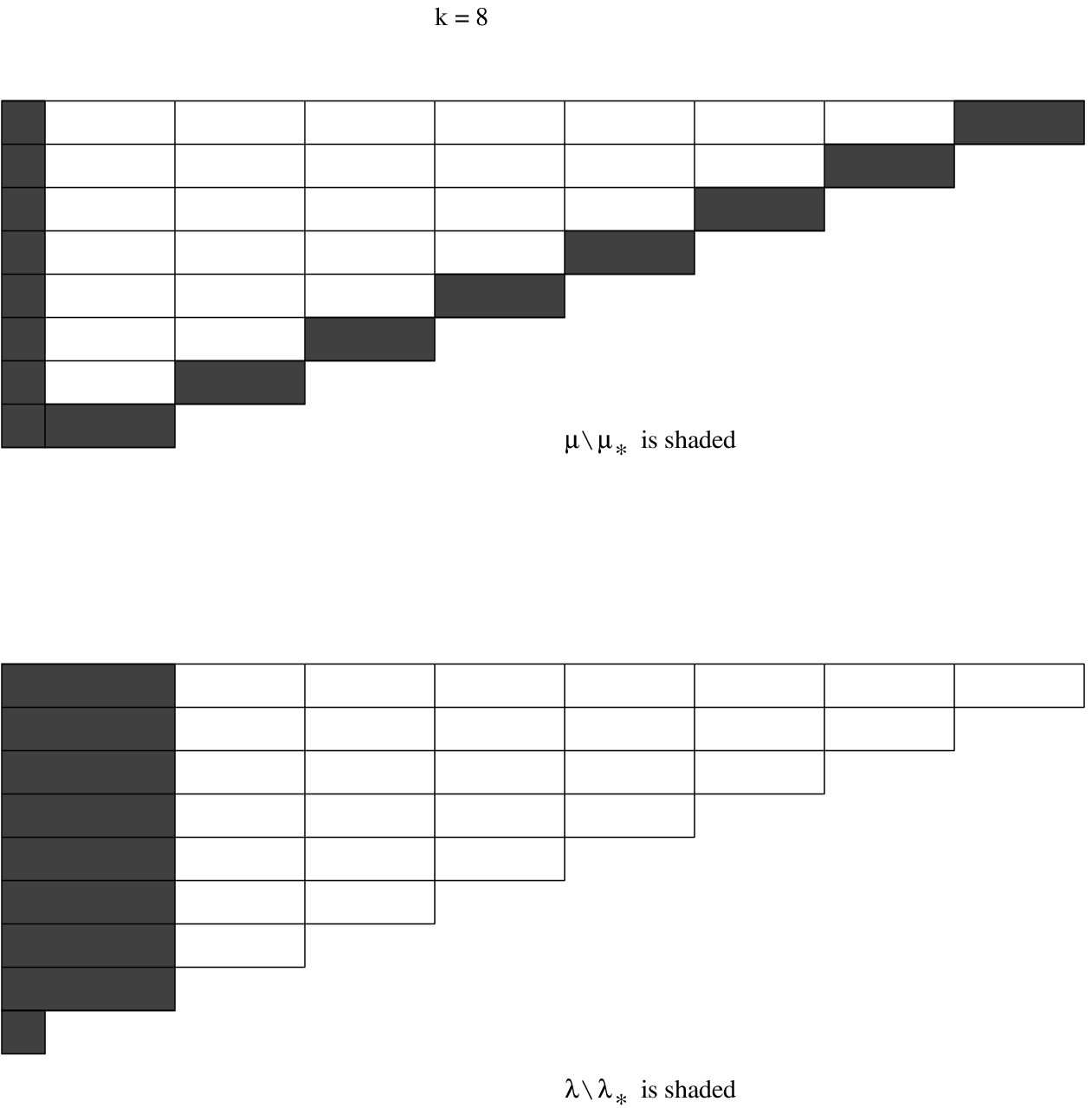}}
\vskip 5pt
\endinsert
\medskip

Denote (and abbreviate)
$$
\HH_{\theta}(\mu)=\prod_{b\in\mu}h_{\mu}(b)=P(b\in\mu)=
P(b\in\mu_*)\cdot P(b\in\mu\backslash \mu_*)\,,
$$
where $P(b\in\mu)$ is the product - of the numbers $h_{\mu}(b)$ - over
the boxes $b\in\mu$, while $P(b\in\mu_*)$ is such a product  - over
the boxes $b\in\mu_*$, etc.

Similarly, let $P'(b \in \lambda)$ denote the product - of the numbers 
$h^{\prime}_{\lambda}(b)$ - over
the boxes $b\in\lambda$, etc., so that
$$
\HH^{\prime}_{\theta}(\lambda)=\prod_{b\in\lambda}h^{\prime}_{\lambda}(b)
=P^{\prime} (b\in\lambda)= P^{\prime}(b\in\lambda_*)\cdot
P^{\prime}(b\in\lambda\backslash \lambda_*)\,.
$$
Thus,
$$
\HH_{\theta}(\mu)
\HH^{\prime}_{\theta}(\lambda)=\left (P(b\in\mu_*)\cdot P^{\prime}
(b\in\lambda_*)\right )\cdot P(b\in\mu\backslash\mu_*)\cdot P
^{\prime}(b\in\lambda\backslash\lambda_*)\,.\tag 2.4.5
$$
We establish below that
$$
\multline
P(b\in\mu_*) \cdot P^{\prime}(b\in\lambda_*)
{\buildrel\hbox{\rm def}\over =}
\prod_{b\in\mu_*} h_{\mu}(b)\cdot \prod_{b\in\lambda_*}h^{\prime}_{\mu}
(b) \\
= \prod_{1\leq i\leq r\leq k-1}\frac {(\mu_i-\mu_r+(r+1-i)\theta)_{\mu_r
-\mu_{r+2}+1}} {\mu_i-\mu_{r+1}+(r+1-i)\theta}\,,
\endmultline
\tag 2.4.6
$$
$$
P(b\in\mu\backslash\mu_*){\buildrel\hbox{\rm def}\over =}
\prod_{b\in\mu\backslash\mu_*} h_{\mu}(b) = \prod^k_{i=1}(\mu_i-\mu_{i+1})
!\cdot \prod^k_{i=1}(\mu_i+\theta(k-i)),
\tag 2.4.7
$$
and that
$$
\multline
P^{\prime}(b\in\lambda\backslash\lambda_*) {\buildrel\hbox{\rm def}\over =}
\prod_{b\in\lambda\backslash\lambda_*} h^{\prime}_{\lambda}(b)\\
= \theta\cdot \prod^k_{i=1} \frac{(\mu_i-\mu_k+(k+1-i)\theta)_{\mu_k-\mu
_{k+2}+1}} {\mu_i-\mu_{k+1}+(k+1-i)\theta} \cdot
\frac{\mu_i-1+(k+2-i)\theta} {\mu_i+(k-i+1)\theta }\  .
\endmultline
\tag 2.4.8
$$
Thus, (2.4.6), (2.4.7) and (2.4.8) clearly imply that
$$
\multline
\HH_{\theta}(\mu)\cdot \HH^{\prime}_{\theta}(\lambda)  \\
=\frac{\theta\cdot\prod_{1\leq i\leq r\leq k}(\mu_i-\mu_r+(r+1-i)\cdot
\theta)_{\mu_r-\mu_{r+2}+1} \cdot \prod^k_{i=1}(\mu_i-\mu_{i+1} )!}
{\prod_{1\leq i\leq r\leq k}(\mu_i-\mu_{r+1}+(r+1-i) \theta)} \times \\
\times \prod^k_{i=1} \frac{ (\mu_i+(k-i)\theta )(\mu_i-1+(k-i+2)
\theta)} {\mu_i+(k-i+1)\theta}\,.
\endmultline
\tag 2.4.9
$$
Substitute $s=r+1$ so that
$$
\multline
\prod_{1\leq i\leq r\leq k}(\mu_i-\mu_{r+1}+(r+1-i)\theta) =\\
= \prod_{1\leq i<s\leq k+1}(\mu_i-\mu_s+(s-i)\theta)  \\
= \prod_{1\leq i<s\leq k}(\mu_i-\mu_s+(s-i)\theta) \cdot \prod^k
_{i=1}(\mu_i-1 +(k+1-i)\theta)\ .
\endmultline \tag 2.4.10
$$
Thus (2.4.9) can be rewritten as $\HH_{\theta}(\mu)\cdot \HH^{\prime}_{\theta}
(\lambda)=A\cdot B$, where
$$
A=\frac {\theta \cdot \prod^k_{i=1}(\mu_i-\mu_{i+1})!\cdot \prod_{1\leq i<s\leq
k+1}((s-i)\theta +\mu_i-\mu_{s-1})_{\mu_{s-1}-\mu_{s+1}+1}}
{\prod_{1\leq i<s\leq k}((s-i)\theta +\mu_i-\mu_s)}
$$
and
$$
B=\prod^k_{i=1} \frac {(\mu_i-1+(k-i+2)\theta)\cdot (\mu_i+(k-i)\theta)}
{(\mu_i+(k-i+1)\theta)\cdot (\mu_i-1+(k-i+1)\theta)}\ .
$$
Note that by (2.1.2), $B=\vk_{\theta}(\mu,\lambda)$, hence
$$
\multline
\frac {\vert\mu\vert !\cdot\theta^{\vert\mu\vert +1}} {\HH_{\theta}(\mu)
\cdot \HH^{\prime}_{\theta}(\lambda)} \vk_{\theta}(\mu,\lambda) =
\frac {\vert\mu\vert !\cdot \theta^{\vert\mu\vert +1}} {A} \\
= \frac {\vert\mu\vert !\cdot\theta^{\vert\mu\vert}\cdot\prod_{1\leq i <s\leq
k} ((s-i)\theta +\mu_i-\mu_s)}
{\prod^k_{i=1} (\mu_i-\mu_{i+1})! \cdot
\prod_{1\leq i<s\leq k+1}((s-i)\theta +\mu_i-\mu_{s-1})_{\mu_{s-1}-
\mu_{s+1}+1}}\,,
\endmultline\tag 2.4.11
$$
and formula (2.4.1) follows.
\bigskip

It remains to verify (2.4.6), (2.4.7) and (2.4.8). \enddemo

\demo{Proof of {\rm(2.4.6)}} Let $1\leq i\leq r\leq k-1$ (hence $\Delta
(i,r+1)\subseteq \mu_*$ and $\Delta(i,r)\subseteq \lambda_*$) and denote
$$
P_{i,r}{\buildrel\hbox{\rm def}\over =} \prod_{b\in\Delta(i,r)}
h^{\prime}_{\lambda}(b)\cdot \prod_{b\in\Delta(i,r+1)}h_{\mu}(b)\ ,
\tag 2.4.12
$$
then it suffices to show that
$$
P_{i,r}= \frac {(\mu_i-\mu_r+(r+1-i)\theta)_{\mu_r-\mu_{r+2}+1}}
{\mu_i-\mu_{r+1}+(r+1-i)\theta}\ .
\tag 2.4.13
$$
So, let $b=(i,j)\in\Delta(i,r)$, then $a_{\mu}(b)=\mu_i-j$ and $\ell_{\mu}
(b)=r-i$, therefore $h_{\theta}(b)=\mu_i-j+(r-i)\theta +1$. Replacing $r$
by $r+1$ we get
$$
\multline
\prod_{b\in\Delta(i,r+1)}h_{\theta}(b) = \prod^{\mu_{r+1}}_{j=\mu_{
r+2}+1} (\mu_i-j+(r+1-i)\theta +1) \\
= (\mu_i-\mu_{r+1}+(r+1-i)\theta +1)_{\mu_{r+1}-\mu_{r+2}}=_{\left ((z+1)_n
= \frac{1}{z}\cdot (z)_{n+1}\right )}\\
= \frac {(\mu_i-\mu_{r+1}+(r+1-i)\theta)_{\mu_{r+1}-\mu_{r+2}+1}}
{\mu_i-\mu_{r+1}+(r+1-i)\theta}\,.
\endmultline
\tag 2.4.14
$$
Similarly, let $b=(i,j)\in\Delta(i,r)$ ($r\leq k-1$), then $h^{\prime}
_{\lambda}(b)=a_{\lambda}(b)+\theta\cdot\ell_{\lambda}(b)+\theta$.
Trivially, $a_{\lambda}(b)=a_{\mu}(b)=\mu_i-j$. Since $r\leq k-1$, hence
$r+1\leq k$, so $1\leq\mu_k\leq\mu_{r+1}$ and $2\leq\mu_{r+1}+1\leq j$.
Therefore $\ell_{\lambda}(b)=\ell_{\mu}(b)=r-i$, so $h^{\prime}_{\theta}
(b)=\mu_i-j+(r+1-i)\theta$. Thus
$$
\prod_{b\in\Delta(i,r)}h^{\prime}_{\lambda}(b)=\prod^{\mu_r}_{j=\mu_{r+1}
+1} (\mu_i-j+(r+1-i)\theta) = (\mu_i-\mu_r+(r+1-i)\theta)_{\mu_r-
\mu_{r+1}}\ .
\tag 2.4.15
$$
Conclude that for $1\leq i\leq r\leq k-1$, by (2.4.14) and (2.4.15),
$$
\multline
\prod_{b\in\Delta(i,r)} h^{\prime}_{\lambda}(b)\cdot\prod_{b\in\Delta
(i,r+1)}h_{\mu}(b) = \\
= \frac {(\mu_i-\mu_r+(r+1-i)\theta)_{\mu_r-\mu_{r+1}}\cdot
(\mu_i-\mu_{r+1}+(r+1-i)\theta)_{\mu_{r+1}-\mu_{r+2}+1}}
{\mu_i-\mu_{r+1}+(r+1-i)\theta} \\
= \frac {(\mu_i-\mu_r+(r+1-i)\theta)_{\mu_r-\mu_{r+2}+1}}
{\mu_i-\mu_{r+1}+(r+1-i)\theta}\ ,
\endmultline
$$
which proves (2.4.6).
\enddemo

\demo{Proof of {\rm(2.4.7)}} Obviously, 
$$
\prod_{b\in\mu\backslash\mu_*}
h_{\mu}(b) = \left [\prod^k_{i=1}\prod_{b\in\Delta(i,i)} h_{\mu}
(b)\right ]\cdot \left [\prod^k_{i=1}h_{\mu}(i,1)\right ]
$$ 
(note that
$(k,1)\notin\Delta (k,k)$, but $(k,1)$ belongs to the left column of $\mu$).
If $b\in\Delta(i,i)$ then $\ell_{\mu}(b)=0$, therefore
$$
\prod_{b\in\Delta(i,i)}h_{\mu}(b) = \prod^{\mu_i}_{j=\mu_{i+1}+1}(\mu_i
-j+1) = (\mu_i-\mu_{i+1})!\ .\tag 2.4.17
$$
Also
$$
h_{\mu}(i,1) = a_{\mu}(i,1)+\theta\ell_{\mu}(i,1)+1 = \mu_i-1 +
(k-i)\theta +1 = \mu_i+(k-i)\theta\ .
$$
Deduce that
$$
\prod_{b\in\mu\backslash\mu_*} h_{\mu}(b) = \prod^k_{i=1}(\mu_i-\mu
_{i+1})! \cdot \prod^k_{i=1}(\mu_i+(k-i)\theta)\ ,\tag 2.4.18
$$
which proves (2.4.7).
\enddemo

\demo{Proof of {\rm(2.4.8)}} Recall that $\Delta(k+1,k+1)=\{(k+1,1)\}$ and
$h^{\prime}_{\lambda}(k+1,1)=\theta$ since $a_{\lambda}(k+1,1)=\ell_{\lambda}
(k+1,1)=0$. Note also that $\Delta(i,k)\sqcup \Delta(i,k+1) = \{(i,j)\in
\lambda\vert 1\leq j\leq \mu_k\}$. Let $1\leq i\leq k$, $1\leq j\leq \mu_k$,
then $h^{\prime}_{\lambda}(i,j)=(\mu_i-j)+\theta(k-i)+\theta$, hence
$$
\prod^{\mu_k}_{j=2} h^{\prime}_{\lambda}(i,j) = \prod^{\mu_k}_{j=2} (\mu_i
-j+(k+1-i)\theta) = (\mu_i-\mu_k+(k+1-i)\theta)_{\mu_k-1}.$$
If $b=(i,1)\in\lambda$ with $1\leq i<k,$ then, 
with respect to $\lambda,$
$$
h^{\prime}_{\lambda}(i,1) = (\mu_i-1)+(k+1-i)\theta+\theta = \mu_i-1 +
(k-i+2)\theta .
$$
Thus
$$
\multline
\prod^{\mu_k}_{j=1} h^{\prime}_{\lambda}(i,j) = (\mu_i-\mu_k+(k+1-i)
\theta)_{\mu_k-1}\cdot (\mu_i-1+(k-i+2)\theta) = \\
= \frac {(\mu_i-\mu_k+(k+1-i)\theta)_{\mu_k-\mu_{k+2}+1}}
{\mu_i-\mu_{k+1}+(k+1-i)\theta} \cdot
\frac {\mu_i-1+(k-i+2)\theta} {\mu_i+(k-i+1)\theta}
\endmultline \tag 2.4.19
$$
(since $\mu_{k+1}=1$, $\mu_{k+2}=0$ and $(u)_n=\frac {(u)_{n+2}}
{(u+n)(u+n+1)}$). Therefore
$$
\multline
\prod_{b\in\lambda\backslash\lambda_*} h^{\prime}_{\lambda}(b) = 
h^{\prime}_{\lambda}(k+1,1)\cdot \prod_{i=1}^k
\prod_{j=1}^{\mu _k} h^{\prime}_{\lambda}(i,j) \\
=\theta\cdot
\prod^k_{i=1} \frac {(\mu_i-\mu_k+(k+1-i)\theta)_{\mu_k-\mu_{k+2}+1}}
{\mu_i-\mu_{k+1}+(k+1-i)\theta} \cdot \frac {\mu_i-1+(k-i+2)\theta}
{\mu_i+(k-i+1)\theta}\ ,
\endmultline
$$
which verifies (2.4.8). \enddemo

This completes the proof of Theorem 2.4.1.

\subhead 2.5. The $z$-measures identity for $(i,j)=(k+1,1)$\endsubhead
The  $z$-measures and $(i,j)=(k+1,1)$ imply:

\proclaim{Theorem 2.5.1}
For the box $(i,j)=(k+1,1),$ the  $z$-measures on the Jack graph
imply the following identity:
$$
\multline
\sum_{\mu _1 \ge \dots \ge \mu _{k+1} =1}
\frac{|\mu |! \cdot \prod _{1 \le i < j \le k}((j-i)\theta +\mu _i-\mu _j)}
{\theta \cdot \prod _{1 \le i \le k} (\mu _i-\mu _{i+1})! \cdot
\prod _{1 \le i < j \le k+1} ((j-i)\theta +
\mu _i-\mu _{j-1})_{\mu _{j-1}-\mu _{j+1}+1}} \\
\times \frac{(z-k\theta )(\bar z-k\theta ) \cdot \prod_{1 \le i \le k}
\left[ (z-(i-1)\theta )_{\mu _i}
(\bar z-(i-1)\theta )_{\mu _i} \right]}
{(\theta ^{-1}z\bar z)_{|\mu |+1}} =1\ ,
\endmultline \tag 2.5.1
$$
where $| \mu |=\mu _1+ \dots + \mu _k\quad (=n-1)\ .$
\endproclaim

\demo{Proof} The proof follows from Theorem 1.8.1., i.e., from the identity
$$
\sum_{\mu} \dim_{\theta}( \mu ) \cdot \varkappa_{\theta} ( \mu, \lambda) \cdot
\varphi_{z\bar z}( \lambda)=1
$$
by arguments similar to those of the previous theorem. Note that here
the function $\varphi_{z\bar z}( \lambda)$ in the $(k+1,1)$ case equals
$$
\varphi_{z\bar z}( \lambda)= \frac
{(z-k\theta )(\bar z-k\theta ) \cdot \prod_{1 \le i \le k} \left[ (z-(i-1)\theta )_{\mu _i}
(\bar z-(i-1)\theta )_{\mu _i} \right]}{(\theta ^{-1}z\bar z)_{|\mu |+1}}\ .
\tag 2.5.2
$$
\enddemo

\subhead 2.6. Some special cases\endsubhead
When $(k+1,1)=(2,1),$
simple  manipulations in the last theorem yield the amazing identity:
$$
{}_3F_2(z+1,\bar z+1, 2; \theta +2, \theta^{-1}z \bar z+2;1)=
\frac{(\theta +1)(z \bar z+\theta )}
{(z-\theta )(\bar z- \theta )}\,. \tag 2.6.1
$$
A long list of cases when the hypergeometric series 
$\;{}_3F_2(a_1,a_2,a_3; b_1,b_2;1)\;$ can be summed in a closed form is
given in the handbook \cite{PBW}, subsection 7.4.4. At a first glance,
(2.6.1) does not appear in that list. However, it is hidden in
formula No.~27, which is as follows:
$$
\multline
{}_3F_2(1,a,b;c,d;1)+ \frac{ab(2+a+b-c-d)}{cd[ab-(c-1)(d-1)]} \cdot
{}_3F_2(2,a+1,b+1;c+1,d+1;1)= \\
\frac {(1-c)(d-1)}{ab-(c-1)(d-1)}\ ,
\endmultline \tag 2.6.2
$$
where $Re(a+b-c-d) < -2$.

To get (2.6.1) from (2.6.2), multiply both sides of (2.6.2) by
$ab-(c-1)(d-1)$ so it cancels from the denominators,
then substitute
$a=z$, $b= \bar z$, $c= \theta +1$, $d= \theta^{-1}z \bar z +1.$
Under this substitution $ab-(c-1)(d-1)=0,$ hence the first summand on
the left vanishes, and the identity (2.6.1) clearly follows.
\medskip
For $k \ge 2,$ no direct proofs
of the identities given by Theorems 2.4.1 and 2.5.1 are known at the
moment.

We now list some special cases of Theorem 2.4.1. Also included
is the $\theta$-Plancherel identity for $(i,j)=(2,2)$, which was
obtained directly from Theorem 1.8.1.

$\bullet$ The box $(i,j)=(2,2)$ implies the identity
$$
\multline
\sum_{r,s=0}^{\infty} \frac{(r+s+3)!\cdot \theta ^{r-s+2}}{[r+2+(s+1)\theta ]
(r+1)!(\theta ^{-1})_{s+1}[r
+1+(s+2)\theta ][r+2\theta ][1+(s+1)\theta ](\theta )_rs!}\times\\
\frac{(r+2\theta )(r+1)}{(r+1+\theta )(r+\theta )}=1\,.
\endmultline \tag 2.6.3
$$

$\bullet$ The box $(i,j)=(2,1)$ implies
$$
\sum_{r \ge 0} \frac{(r+1) \cdot\theta ^{r+1}}{(\theta   )_{r+2}}
=\sum_{r \ge 0} \frac{(r+1)! \cdot\theta  ^{r+1}}{r! \cdot(\theta   )_
{r+2}} =1. \tag 2.6.4
$$
It is easy to give (2.6.4) a direct  proof.

$\bullet$ The box $(i,j)=(3,1)$ implies
$$
\sum_{s \ge r \ge 0} \frac{(s+r+2)! \cdot\theta   ^{s+r+2} \cdot (\theta
 +s-r)}
{r! \cdot (s-r)! \cdot (\theta  )_{r+2}\cdot (\theta  )_{s+1}(2\theta
 +s-r)_{r+2}}=1. \tag 2.6.5
$$

$\bullet$ The box $(i,j)=(4,1)$ implies that
$\sum_{u \ge s \ge r \ge 0}f_4(r,s,u;\theta  )=1$, where
$$
\multline
f_4(r,s,u;\theta   )= \\
\frac{(r+s+u+3)! \cdot \theta   ^{r+s+u+3}(\theta   +s-r)(\theta   +u-s)
(2\theta   +u-r)}
{r! \cdot (s-r)! \cdot (u-s)! (t)_{r+2}(\theta   )_{s+1}(\theta
  )_{u-r+1}} \times \\
\frac{1}{(2\theta   +s-r)_{r+2}(2\theta   +u-s)_{s+1}(3\theta
  +u-r)_{r+2}}\ . \endmultline \tag 2.6.7
$$

$\bullet$ The box  $(i,j)=(5,1)$ implies that
$\sum_{v \ge u \ge s \ge r \ge 0}f_5(r,s,u,v;t)=1$, where
$$
\multline
f_5(r,s,u,v;\theta   )= \\
\frac{(r+s+u+v+4)! \cdot \theta   ^{r+s+u+v+4}(\theta   +s-r)(t+u-s)(2\theta
  +u-r)}
{r! \cdot (s-r)! \cdot (u-s)! (\theta   )_{r+2}(\theta   )_{s+1}( \theta
 )_{u-r+1}} \times \\
\frac{1}{(2\theta   +s-r)_{r+2}(2\theta   +u-s)_{s+1}(3\theta
  +u-r)_{r+2}} \times \\
\frac{(\theta   +v-u)(2 \theta  +v-s)(3\theta   +v-r)}
{(v-u)!( \theta  )_{v-s+1}(2\theta   +v-u)_{u-r+1}(3\theta   +v-s)_{s+1}
(4\theta   +v-r)_{r+2}}\ .
\endmultline \tag 2.6.8
$$

\subhead 2.7. The case $\theta=1$ \endsubhead
In this case the measure $M_\infty$ becomes the Plancherel measure
for the Young graph $\BbbY$, and the identity (2.4.1), which corresponds to
the box $(k+1,1)$, takes a simpler form given in (0.6). On the other hand,  
the identity for the Plancherel measure on $\BbbY$ and a general box
$(k+1,l+1)$ was given in \cite{Re}, see (0.2) above. The
specialization of (0.2) for $l=1$ takes the form (0.7). Here we
verify directly that the left hand sides of these two identities,
(0.6) and (0.7), coincide, thus proving they indeed are the same identity.
\bigskip

\proclaim{Proposition 2.7.1} Let
$$
f_k(p_1 , \dots , p_k )=
\frac{(p_1+ \dots + p_k -k(k-1)/2)! \cdot
\prod_{1 \le i < j \le k}(p_i-p_j)^2}
{ \prod_{i=1}^k[(p_i-1)!(p_i+1)!] } \tag 2.7.1
$$
so that (0.6) is the identity
$$
\sum_{p_1 > \dots > p_k \ge 1} f_k(p_1 , \dots , p_k )=1. 
$$
Similarly, let $\mu_{k+1}=1$ and
$$
g_k(\mu_1, \dots ,\mu _k)=
\frac{(\mu_1+ \dots +\mu _k)! \prod_{1 \le i < j \le k}(\mu_i-\mu_j+j-i)}
{\prod_{i=1}^k(\mu_i - \mu_{i+1})! \cdot
\prod_{1 \le i < j \le k+1}(\mu_i-\mu_{j-1}+j-i)_{\mu_{j-1}-\mu_{j+1}+1}}
\tag 2.7.2 $$
so that (0.7) is the identity
$$
\sum_{\mu_1 \ge \dots \ge \mu_{k+1}=1}g_k(\mu_1, \dots ,\mu _k)=1.
$$
Finally, let $p_j=\mu_j+k-j, \quad 1 \le j \le k+1$
(so $p_{k+1}=0$), $p_r=0$ if $r \ge k+1.$

Then the expressions (2.7.1) and (2.7.2) are equal. Hence
(0.6) and (0.7) indeed are the same identity.
\endproclaim

\demo{Proof} Define $f_k^*(p)$ via
$$
g_k(\mu_1, \dots ,\mu _k)=
g_k(p_1-(k-1), \dots , p_k)=f_k^*(p_1, \dots , p_k).\tag 2.7.3
$$

The statement of the above proposition clearly follows from

\proclaim{Claim 1} $f_k(p_1 , \dots , p_k )=f_k^*(p_1, \dots , p_k).$
\endproclaim

\demo{Proof of Claim 1} First calculate $f_k^*(p_1 , \dots , p_k)$.
We have
$$
\gather
\mu_1+ \dots +\mu_k=p_1+ \dots +p_k -k(k-1)/2, \qquad
\mu_i-\mu_j+j-i=p_i-p_j,\\
\mu_i-\mu_{i+1}=p_i-p_{i+1}-1
\endgather
$$ 
and
$$
\multline
\prod_{1 \le i < j \le k+1}(\mu_i-\mu_{j-1}+j-i)_{\mu_{j-1}-\mu_{j+1}+1}\\
= \prod_{1 \le i < j \le k}(p_i-p_{j-1}+1)_{p_{j-1}-p_{j+1}-1}
\prod_{i=1}^k(p_i-p_k+1)_{p_k+1}\,.
\endmultline
\tag 2.7.4
$$
Thus,
$$
\multline
f_k^*(p)
=\frac{(p_1+ \dots +p_k-k(k-1)/2)! \prod_{1 \le i < j \le k}(p_i-p_j)}
{\prod_{i=1}^k(p_i-p_{i+1}-1)! \prod_{1 \le i < j \le k}
(p_i-p_{j-1}+1)_{p_{j-1}-p_{j+1}-1}}\\
\times\frac1{\prod_{i=1}^k(p_i-p_k+1)_{p_k+1}}\,. 
\endmultline
\tag 2.7.5
$$
Comparing $f_k(p)$ with $f_k^*(p)$ and cancelling
$(p_1+ \dots + p_k-k(k-1)/2)! \prod_{1 \le i<j \le k}(p_i-p_j)$,
we see that Claim 1 is equivalent to \enddemo

\proclaim{Claim 2}
$$ \multline
\prod_{1 \le i < j \le k}[(p_i-p_j) \cdot 
(p_i-p_{j-1}+1)_{p_{j-1}-p_{j+1}-1}]
\prod_{i=1}^{k}(p_i-p_{i+1}-1)!(p_i-p_k+1)_{p_k+1}= \\
=\prod_{i=1}^{k}(p_i-1)!(p_i+1)! \endmultline \tag 2.7.6
$$
\endproclaim

Note that in the second factor on the left, when $i=k$, the corresponding
term in that product becomes $(p_k-1)!(p_k+1)!$. This same term also appears
on the right hand side --- again when $i=k$.  Cancelling it on both sides we
see that Claim 2 is equivalent to 

\proclaim{Claim 3}
$$ \multline
\prod_{1 \le i < j \le k}[(p_i-p_j) \cdot (p_i-p_{j-1}+1)_{p_{j-1}-p_{j+1}-1}]
\prod_{i=1}^{k-1}(p_i-p_{i+1}-1)!(p_i-p_k+1)_{p_k+1}= \\
=\prod_{i=1}^{k-1}(p_i-1)!(p_i+1)! \endmultline \tag 2.7.7
$$
\endproclaim

On the left--hand side  of (2.7.7) the terms involving the index $i=1$ are
$$
R_1=\left [\prod_{j=2}^k(p_1-p_j)(p_1-p_{j-1}+1)_{p_{j-1}-p_{j+1}-1}\right]
 \cdot (p_1-p_2-1)!(p_1-p_k+1)_{p_k+1}.\tag 2.7.8
$$
By induction on $k$, the proof of Claim 3 will follow once we prove 

\proclaim{Claim 4}
$$
R_1=(p_1-1)!(p_1+1)!
$$
\endproclaim

We need 

\proclaim{Claim 5}  Let $2 \le r \le k-1$, then
$$
(p_1-p_2-1)!\cdot \left [\prod_{j=2}^r(p_1-p_j)
(p_1-p_{j-1}+1)_{p_{j-1}-p_{j+1}-1}\right ]
=(p_1-p_r)!(p_1-p_{r+1}-1)!\tag 2.7.9
$$
\endproclaim

\demo{Proof of Claim 5} Induction on $r$, where $2 \le r \le k-1$.

$r=2: \quad (p_1-p_2)(p_1-p_1+1)_{p_1-p_3-1}
(p_1-p_2-1)!=(p_1-p_2)!(p_1-p_3-1)!$
\medskip
$r \Rightarrow r+1$:
$$ \gather
\prod_{j=2}^{r+1} \dots = (\prod_{j=2}^ r \dots) \times
(p_1-p_{r+1})(p_1-p_r+1)_{p_r-p_{r+2}-1}  \quad \text{by induction} \\
=(p_1-p_r)!(p_1-p_{r+1}-1)!\times
(p_1-p_{r+1})(p_1-p_r+1)_{p_r-p_{r+2}-1}\\
=(p_1-p_{r+1})!(p_1-p_{r+2}-1)!\tag2.7.10
\endgather
$$
which proves Claim 5.
\enddemo

We can now prove Claim 4, thus completing the proof of proposition 2.7.1.

Rearrange terms in $R_1$ and apply Claim 5 with $r=k-1$:
$$
\multline
R_1
=\{\prod_{j=2}^{k-1}[(p_1-p_j)(p_1-p_{j-1}+1)_{p_{j-1}-p_{j+1}-1}]
(p_1-p_2-1)!\}\\
\times (p_1-p_k)(p_1-p_{k-1}+1)_{p_{k-1}-1}(p_1-p_k+1)_{p_k+1}\\
=\{(p_1-p_{k-1})!(p_1-p_k-1)!\}(p_1-p_k)(p_1-p_{k-1}+1)_{p_{k-1}-1}
(p_1-p_k+1)_{p_k+1}\\
=[(p_1-p_{k-1})!(p_1-p_{k-1}+1)_{p_{k-1}-1}][(p_1-p_k-1)!(p_1-p_k)
(p_1-p_k+1)_{p_k+1}]\\
=(p_1-1)![(p_1-p_k)!(p_1-p_k+1)_{p_k+1}]\\
=(p_1-1)!(p_1+1)!
\endmultline
\tag 2.7.11
$$

This completes the proof of Proposition 2.7.1. \qed
\enddemo

\subhead 2.8. Remark on \cite{MMW} \endsubhead
Recall that the space $\Tab$ can be viewed as the projective limit
space $\varprojlim\Tab_n$.  Let $M$ be the Plancherel measure
($\theta=1$) on $\Tab$, let $M_n$ be the pushforward of $M$ under the
projection $\Tab\to\Tab_n$, and let $M'_n$ be the uniform probability
measure on $\Tab_n$. Note that the measures $M_n$ and $M'_n$ are
quite different. Indeed, 
$$
M_n(\Tab(\lambda))=\operatorname{const}_1\,\cdot\, (\dim\lambda)^2, \qquad
M'_n(\Tab(\lambda))=\operatorname{const}_2\,\cdot\, \dim\lambda.
$$
Nevertheless, as is shown by McKay, Morse, and Wilf, 
$$
\lim_{n\to\infty}M'_n=M, \tag2.8.1
$$
see \cite{MMW, Section 4}. From this fact McKay, Morse, and Wilf
deduce the following result. Let $T_n\in\Tab_n$ be the random (with
respect to $M'_n$) finite tableau and let $p(n;i,j,k)$ stand for the
probability that $T_n$ has the entry $k$ at the box $(i,j)$. Then
the limit $p_{i,j}(k)=\lim_{n\to\infty}p(n;i,j,k)$ exists and coincides
with the probability $\PP_M(T(i,j)=k)$ (which can be calculated as
is shown in \cite{Re}). 

An explanation of the equality (2.8.1) can be extracted from \cite{V}.

\head \S3 The Kingman graph case\endhead

\subhead 3.1. Kingman edge multiplicities\endsubhead
As the Jack parameter $\theta >0$ goes to 0 (denoted $ \theta\downarrow 0$),
 the Jack symmetric functions
$P_{\mu}$ degenerate to the monomial symmetric functions $m_{\mu}$ \cite{Ma}.
The simplest case of Pieri's formula for the monomial symmetric functions
has the form
$$
m_{\mu} m_{(1)} =\sum_{\lambda:\lambda\searrow\mu}\vk_0(\mu,
\lambda) m_{\lambda}\ ,\tag 3.1.1
$$
where $\vk_0(\mu,\lambda)$ is a strictly positive integer. Specifically,
if $k$ stands for the length of the row in $\lambda$ containing the box
$\lambda\backslash\mu$, then $\vk_0(\mu,\lambda)$ is the multiplicity
of $k$ in $\lambda$.

One can verify that
$$
\vk_0(\mu,\lambda)=\lim_{\theta\downarrow 0} \vk_{\theta}(\mu,\lambda)\ ,
\quad \mu\nearrow\lambda\ .\tag 3.1.2
$$
We take the numbers $\vk_0(\mu,\lambda)$ as edge multiplicities. The graph
$\Bbb Y$ together with these edge multiplicities is called the {\it Kingman
graph,\/} see \cite{Ke1, KOO, BO3}. The dimension function of the
Kingman graph will be denoted as $\dim_0\mu$. It is given by a simple
formula 
$$
\dim_0\mu =\frac{\vert\mu\vert !}{\mu_1!\mu_2!\dots \mu_{\ell}!}\,,
\tag 3.1.3
$$
where $\ell=\ell(\mu)$ is the number of nonzero rows in $\mu$.

Again, one can verify that
$$
\dim_0\mu =\lim_{\theta\downarrow 0} \dim_{\theta}\mu\ .
$$

\subhead 3.2. The $t$-measures\endsubhead
We fix a parameter $t>0$.

\proclaim{Theorem 3.2.1} For any $t>0$ there exists a $\vk_0$-central measure
$M_t$ such that the corresponding $\vk_0$-harmonic function has the form
$$
\psi_t(\lambda)=\frac{(\lambda_1-1)! \dots (\lambda_{\ell}-1)!}
{r_1(\lambda)!r_2(\lambda)!\dots}\
\cdot \frac{t^{\ell}}{(t)_n}\ ,\tag 3.2.1
$$
where $\ell=\ell(\lambda)$ and $r_k(\lambda)$ is the multiplicity of $k$ in
$\lambda$.
\endproclaim

\demo{Proof} See \cite{BO3}. \qed
\enddemo

\subhead 3.3. Reachability for $t$-measures\endsubhead

\proclaim{Lemma 3.3.1} Let $M=M_t$. Then all the boxes are completely
reachable.
\endproclaim

\demo{Proof} It is readily verified that Lemma 2.3.3 holds for the measure
$M$. Then we apply the same argument as in the proof of Theorem
2.3.1.\qed
\enddemo

\proclaim{Theorem 3.3.2} Let $k=0,1,\dots$, $l=1,2,\dots$, and $t > 0$.
Then
$$
\multline
\frac{1}{k!} \sum_{r_1,\dots,r_l,s_1, \dots,s_k \geq 0}
\frac{(s_1+ \dots +s_k+r_1+2r_2+ \dots +lr_l+kl+k+l)!}
{(s_1+l+1) \dots (s_k+l+1)1^{r_1}2^{r_2} \dots l^{r_l}r_1! \dots r_l!} \\
\times \frac{t^{k+r_1+ \dots +r_l+1}}
{(t)_{s_1+ \dots +s_k+r_1+2r_2+ \dots +lr_l+kl+k+l+1}}
=1. \endmultline \tag 3.3.1
$$
\endproclaim

We prove this in two ways: the first proof is based on Theorem
1.8.1, while the second one is a direct derivation. A crucial step in
that direct proof was shown to us by S.~Milne. \medskip

\subhead 3.4. First proof \endsubhead
Fix the box $(k+1,l+1)$ and apply Theorem 1.8.1 (the complete
reachability is ensured by Lemma 3.3.1).
Since the box $(k+1,l+1)$ can be added to $\mu$ (yielding $\lambda$),
we have 
$$
\gather
\mu_1 \ge \dots \ge\mu_k \ge l+1, \quad \mu_{k+1}=l, \quad
r_l(\mu)\ge1,\\
r_l(\lambda)=r_l(\mu)-1, \quad 
r_{l+1}(\lambda)=r_{l+1}(\mu)+1,\\ 
r_j(\lambda)=r_j(\mu) \quad \text{if $j \neq l,l+1$,} \\
\ell(\lambda)=k+r_1(\mu)+\dots+r_l(\mu). \\
\endgather
$$

Let us abbreviate
$$
\gather
r_1=r_1(\mu), \, r_2=r_2(\mu),\, \dots,\\
\ell=\ell(\lambda)=k+r_1+\dots+r_l
\endgather
$$
and remark that 
$$ 
\varkappa(\mu,\lambda)=r_{l+1}(\lambda)=r_{l+1}+1.
$$
In this notation, the identity 
$$
\sum_{\mu} \dim \mu \cdot \varkappa ( \mu, \lambda )
\cdot \psi_t(\lambda) = 1 
$$
provided by Theorem 1.8.1 becomes
$$
\sum_{r_1, \dots,r_{l-1} \geq 0} \sum_{r_l \geq 1}
\sum_{\mu_1 \geq \dots \geq \mu_k \geq l+1}
\frac{|\mu|!}{\mu_1! \mu_2! \dots}(r_{l+1}+1)
\frac{(\lambda_1-1)!(\lambda_2-1)! \dots}{r_1(\lambda)!r_2(\lambda)!\dots}
\frac{t^\ell}{(t)_{|\mu|+1}}
=1. \tag 3.4.1
$$ \bigskip

Now $\lambda _{k+1} = \mu_{k+1}+1=l+1$ and $\lambda_j=\mu_j$
if $j \neq k+1$. Therefore
$$
\frac{(\lambda_1-1)!(\lambda_2-1)! \dots}{\mu_1! \mu_2! \dots}=
\frac{1}{\mu_1 \dots \mu_k 1^{r_1}2^{r_2} \dots(l-1)^{r_{l-1}}l^{r_l-1}}\,.
\tag 3.4.2
$$
Also,
$$
\frac{r_{l+1}+1}{r_1(\lambda)!r_2(\lambda)! \dots} =
\frac{1}{r_1! \dots r_{l-1}!(r_l-1)!r_{l+1}!r_{l+1}! \dots}\,. \tag 3.4.3
$$
Substituting (3.4.2) and (3.4.3) into (3.4.1) and replacing $r_l-1$ by
$r_l$, we obtain
$$
\sum_{r_1, \dots,r_l \geq 0} \sum_{\mu_1 \geq \dots \geq \mu_k \geq l+1}
\frac{|\mu|!}{\mu_1 \dots \mu_k 1^{r_1} \dots l^{r_l}r_1! \dots r_l!}
\frac{1}{r_{l+1}!r_{l+2}! \dots}
\frac{t^\ell}{(t)_{|\mu|+1}} =1\,, \tag 3.4.4
$$
where 
$$
|\mu|=\mu_1+ \dots +\mu_k+r_1+2r_2+ \dots +lr_l+l
$$ 
and
$$
\ell=\ell(\lambda )  =k+r_1+ \dots + r_l+1
$$
(because we have changed the definition of $r_l$).

Note that $r_{l+1}+r_{l+2}+ \dots = k$. Multiply and divide (3.4.4) by $k!$
and note that  $\dfrac{k!}{r_{l+1}!r_{l+2}!\dots}$ equals the number
of permutations of $\mu_1, \dots, \mu_k$. We can therefore
cancel
the factor $\dfrac{k!}{r_{l+1}!r_{l+2}!\dots}$ by replacing 
`$\mu _1\ge \cdots \ge \mu _k $' by 
`$\mu _1, \cdots , \mu _k $'  under the summation sign,
 and (3.4.4) becomes
$$
\frac{1}{k!} \sum_{r_1, \dots ,r_l \geq 0}
\sum_{\mu_1, \dots , \mu_k \geq l+1}
\frac{|\mu|!}{\mu_1 \dots \mu_k 1^{r_1} \dots l^{r_l}r_1! \dots r_l!}
\frac{t^\ell}{(t)_{|\mu|+1}}=1. \tag 3.4.5
$$
Finally, replace $\mu_j$ by $s_j+l+1,\  s_j \geq 0$ and observe that
$$
|\mu|= s_1 + \dots +s_k+r_1+2r_2+\dots+lr_l+kl+k+l.
$$
The first proof of Proposition 3.3.2 clearly follows. \hfill \qed

\subhead 3.5. Second proof (direct)\endsubhead
We first transform the left
hand side of (3.3.1).

\proclaim{Lemma 3.5.1} Let $L$ denote the left
hand side of (3.3.1). We have
$$
\multline
L=\frac{t^{k+1}}{k!} \int_0^1
(1-v)^{t-1} \cdot v^l \cdot \exp \left[t \left(v+\frac{v^2}2+ \dots
+\frac{v^l}l \right) \right] \times \\
\times \left[-\ln(1-v)- \left(v+\frac{v^2}2+ \dots +\frac{v^l}l
\right) \right]^k\,dv.
\endmultline \tag 3.5.1
$$
\endproclaim

\proclaim{Corollary 3.5.2} The identity $L=1$, i.e, the identity {\rm(3.3.1)},
is equivalent to the following integral identity:
$$
\multline
\int_0^1
(1-v)^{t-1} \cdot v^l \cdot \exp \left[t \left(v+\frac{v^2}2+ \dots
+\frac{v^l}l \right) \right] \times  \\
\times \left[-\ln(1-v)- \left(v+\frac{v^2}2+ \dots +\frac{v^l}l
\right) \right]^k\,dv = \frac{k!}{t^{k+1}}\,.
\endmultline \tag 3.5.2
$$
\endproclaim

\demo{Proof of Lemma 3.5.1} Note first that the Taylor
expansion for $\ln(1-v)$ implies that
$$
\sum_{\mu=l+1}^{\infty} \dfrac{v^{\mu}}{\mu} =
-\ln(1-v)- \left(v+\frac{v^2}2+ \dots +\frac{v^l}l \right). \tag 3.5.3
$$

Return now to (3.3.1) with $\mu_i=s_i+l+1$ and write
$n=n(\mu,r)=\mu_1+ \dots +\mu_k+r_1+2r_2+ \dots +lr_l+l$. Then
$$
L=\frac 1{k!} \sum_{r_1, \dots ,r_l = 0}^ \infty \,
\sum_{\mu_1, \dots, \mu_k = l+1}^\infty
\frac {n!}{(t)_{n+1}} \cdot \frac {t^{k+1}  \cdot
t^{r_1 + \dots + r_l}}{ \mu_1 \dots \mu_k \cdot
1^{r_1} \cdot r_1! \cdot 2^{r_2} \cdot r_2! \dots
l^{r_l} \cdot r_l!}\,. \tag 3.5.4
$$

By Euler's beta integral,
$$
\frac{n!}{(t)_{n+1}} = \int_0^1(1-v)^{t-1} \cdot v^n \,dv,\tag 3.5.5
$$
hence
$$
L=\frac{t^{k+1}}{k!} \int_0^1 \, (1-v)^{t-1} \cdot v^l \cdot
\left[ \sum_{r_1=0}^ \infty \frac{t^{r_1}v^{r_1}}
{1^{r_1}r_1!} \right]^{} \dots \left[
\sum_{r_l=0}^ \infty\frac{t^{r_l}v^{l r_l}}{l^{r_l}r_l!} \right]^{} \cdot
\left[ \sum_{\mu=l+1}^ \infty \frac{v^{\mu}}{\mu} \right]^k \,dv.
\tag 3.5.6
$$

Now $$ \sum_{r_j=0}^ \infty \frac{t^{r_j}v^{j r_j}}{j^{r_j}r_j!}=\exp \left[
\dfrac{t \cdot v^j}j \right],$$ hence by (3.5.3),

$$
\multline
L= \frac{t^{k+1}}{k!} \int_0^1 \,(1-v)^{t-1} \cdot v^l \cdot
\exp \left[ t \left( v+ \frac {v^2}2 + \dots + \frac{v^l}l
\right) \right] \times \\
\times \left[ -\ln (1-v) - \left( v+ \dots + \frac{v^l}l
\right) \right]^k \,dv.
\endmultline \tag 3.5.7
$$
\enddemo
This proves (3.5.1).

We proceed to the proof of (3.5.2), which, by Corollary 3.5.2, is
equivalent to the initial identity. 

Denote $y=v+v/2+ \dots + v^l/l$, then
$ \dfrac{dy}{dv}= \dfrac{1-v^l}{1-v}\,$. Therefore
$$
\frac d{dv} \exp (t \cdot y)= \exp (t \cdot y) \cdot t \cdot
\frac{1-v^l}{1-v}\,. \tag 3.5.8
$$
It follows that
$$
\frac{d}{dv} \left( (1-v)^t \cdot \exp (t \cdot y) \right)=
-t \cdot (1-v)^{t-1} \cdot v^l \cdot \exp (t \cdot y). \tag 3.5.9
$$
Since
$$
\frac d{dv}[- \ln (1-v) - y]= \frac
{v^l}{1-v}\,,
$$
we get
$$
\frac{d}{dv} [- \ln (1-v) - y]^k = k \cdot [- \ln (1-v) - y]^{k-1}
\cdot \frac{v^l}{1-v}\,.  \tag 3.5.10
$$
Thus, for $k \ge 0$
$$
\multline
\frac{d}{dv} \left( (1-v)^t \cdot \exp (t \cdot y) \cdot
[- \ln (1-v) - y]^k \right) = \\
-t \cdot (1-v)^{t-1} \cdot v^l \cdot \exp (t \cdot y) \cdot
[- \ln (1-v) - y]^k + \\
 k \cdot (1-v)^{t-1} \cdot v^l \cdot \exp (t \cdot y) \cdot
[- \ln (1-v) - y]^{k-1}.
\endmultline \tag 3.5.11
$$
Since $t>0$ by the assumption, we have, for any $j$,
$$ 
\lim_{v\to 1} (1-v)^t \cdot (\ln (1-v))^j = 0.
$$ 
Thus,
$$
\multline
\left(1-v)^t \cdot \exp (t \cdot y) \cdot [- \ln (1-v) - y]^k
\right |_0^1 =
\cases  -1, & k=0 \\
0 - 0 = 0, & k > 0  \endcases \qquad \qquad \endmultline \tag 3.5.12
$$

Denote by $a_k$ the integral in the left--hand side of (3.5.2). When
$k=0$, the integration of (3.5.9) easily implies that 
$$
a_0 =\int_0^1 (1-v)^{t-1} \cdot v^l \cdot \exp (t \cdot y) \,dv
= \frac{1}{t}\,, \tag 3.5.13
$$
as required.

For arbitrary $k=1,2,\dots$, integrating the expression (3.5.11)
and using (3.5.12) we get
$$
\multline
a_k=\int_0^1 (1-v)^{t-1} \cdot v^l \cdot \exp (t \cdot y) \cdot
[- \ln (1-v) - y]^k \, dv=  \\
\frac{k}{t} \int_0^1 (1-v)^{t-1} \cdot v^l \cdot \exp (t \cdot y) \cdot
[- \ln (1-v) - y]^{k-1} \,dv,
\endmultline \tag 3.5.14
$$
i.e,
$$ 
a_k= \frac{k}{t} \cdot a_{k-1}\,.
$$
Thus, by induction,
$a_k= \dfrac{k!}{t^{k+1}}\,$. \qed

\bigskip
\bigskip

{\smc Grigori Olshanski:} Dobrushin Mathematics Laboratory, Institute for
Information Transmission Problems, Bolshoy Karetny 19, Moscow 101447,
GSP-4, Russia.

E-mail address: {\tt olsh\@iitp.ru, olsh\@online.ru}

\medskip

{\smc Amitai Regev:} Department of Mathematics, The Weizmann Institute of
Science, Rehovot 76100, Israel.

E-mail address: {\tt regev\@wisdom.weizmann.ac.il}

\enddocument

\end